\documentclass[12pt,a4paper]{amsart}
\usepackage[utf8]{inputenc}	
\usepackage[T1]{fontenc}
\usepackage[in,headings]{fullpage}
\usepackage{amsbsy, amscd, amsfonts, amsmath, amsrefs, amssymb, amsthm}
\usepackage{graphicx}
\usepackage{float}
\usepackage{color}   
\usepackage[colorlinks=true,linkcolor=blue,citecolor=blue,urlcolor=blue]{hyperref}

\newtheorem{theorem}{Theorem}

\newtheorem{corollary}[theorem]{Corollary}
\newtheorem{proposition}[theorem]{Proposition}
\newtheorem{lemma}[theorem]{Lemma}

\theoremstyle{definition}
\newtheorem{definition}[theorem]{Definition}
\newtheorem{remark}[theorem]{Remark}

\theoremstyle{remark}

\newcommand{\C}{\mathbf{C}}

\newcommand{\R}{\mathbf{R}}

\renewcommand{\Re}{\mathop{\mathrm{Re}}\nolimits}
\renewcommand{\Im}{\mathop{\mathrm{Im}}\nolimits}
\newcommand{\Rzeta}{\mathop{\mathcal R }\nolimits}

\newfont{\cmbsy}{cmbsy10}
\newfont{\cmmib}{cmmib10}
\newcommand{\Orden}{\mathop{\hbox{\cmbsy O}}\nolimits}
\newcommand{\orden}{\mathop{\hbox{\cmmib o}}\nolimits}

\usepackage{pifont}
\newcommand{\checked}{\text{\ding{51}\ \ }}%
\renewcommand{\checked}{}%

\DeclareMathOperator*{\per}{P} 

\begin{document}

\title[Siegel results about zeros of $\Rzeta$]
{On Siegel results about the zeros of the auxiliary function of Riemann.}
\author[Arias de Reyna]{J. Arias de Reyna}
\address{%
Universidad de Sevilla \\ 
Facultad de Matem\'aticas \\ 
c/Tarfia, sn \\ 
41012-Sevilla \\ 
Spain.} 




\email{arias@us.es, ariasdereyna1947@gmail.com}


\begin{abstract}
We state and give complete proof of the results of Siegel about the zeros of the auxiliary function of Riemann $\Rzeta(s)$. We point out the importance of the determination of the limit to the left of the zeros of $\Rzeta(s)$ with positive imaginary part, obtaining the term $-\sqrt{T/2\pi}P(\sqrt{T/2\pi})$ that would explain the periodic behaviour observed with the statistical study of the zeros of $\mathop{\mathcal R}(s)$.

We precise also the connection of the position on the zeros of $\Rzeta(s)$ with the zeros of $\zeta(s)$ in the critical line.
\end{abstract}

\maketitle

\section{Introduction} 
In his paper \cite{R} Riemann asserts to have proved that all but an infinitesimal proportion of zeros of zeta are in the critical line. This is repeated in a letter to Weierstrass \cite{R2}*{p.~823--825}, where he says that the difficult proof depends on a new development of $\Xi(t)$ that he has not simplified sufficiently to communicate. The Riemann-Siegel expansion for $\zeta(s)$ were recovered by Siegel from Riemann's papers. Siegel in his publication \cite{Siegel} concluded that in Riemann's Nachlass there is no approach to the proof of his assertion on the real zeros of $\Xi(t)$. However, Siegel connected the zeros of the new function $\Rzeta(s)$, found in Riemann's papers, with the zeros of $\zeta(s)$. The main objective of this paper is to expose and give a complete proof of the results of Siegel on the zeros of $\Rzeta(s)$.  

Siegel \cite{Siegel} proved an asymptotic expression \eqref{E:withUT3} valid for $1-\sigma\ge t^a$ with $a=\frac37$. He claims that we may prove this with $a=\varepsilon$ for any $\varepsilon>0$. In \cite{A98}*{Th.~11} we proved it for $1-\sigma\ge t^{2/5}\log t$. In this paper, it will be desirable to have proved the claim with $a<\frac13$ so that the term in $|\sigma|^{3/2}=\orden(t^{1/2})$.  For this reason we will keep $a$ undetermined, to see the possible consequences of the hypothesis of being $a <\frac13$. On the other hand, we may assume in all our results that $a=\frac{3}{7}$ as done by Siegel (see also \cite{A193}*{Th.~9}).

The paper contains many computations that Siegel does not detail. First, we state three main Theorems in Section \ref{S:main}. We postpone its proof, and in Section \ref{S:results} we give the main applications to the zeros of $\Rzeta(s)$. We note that in \cite{A185} we have given a better result about the number of zeros than the one obtained by Siegel methods (see Remark \ref{R:improve}). Our corollary \ref{C:leftzeros} slightly improves on Siegel's results. This corollary implies that $\Rzeta(s)$ contains at least $>cT$ zeros to the left of the critical line, with $c=35/198\pi$. This constant is a slight improvement over the one in Siegel, in any case it is a very poor result, as it only shows that the number of zeros  of $\zeta(s)$ on the critical line $N_0(T)$ is greater than $cT$. 
Our computation of the zeros of $\Rzeta(s)$ \cite{A172} makes plausible that $\frac23$ of the zeros of $\Rzeta(s)$ are on the left  of the critical line, this would imply that $2/3$ of the zeros of $\zeta(s)$ will be simple zeros on the critical line.

In Section \ref{S:Th3} we prove Theorem \ref{T:logrzeta}. It is in Siegel only implicitly.
In Section \ref{S:Th2} we prove Theorem \ref{T:firstLittlewood}. We make a modification of Siegel's reasoning. He uses a function $g(s)=\pi^{-(s+1)/2}e^{-\pi i s/4}\Gamma(\frac{s+1}{2})\Rzeta(s)$, while we use $F(s)=s\pi^{-s/2}\Gamma(s/2)\Rzeta(s)$. 
The main advantage is that our function is entire, while $g(s)$ has poles and zeros at the negative real axis, making the Siegel reasoning difficult.  We postpone for section \ref{S:comp} the proof of propositions  \ref{P:logF} and \ref{P:intvert} that are mere computations.

\subsection{Notations} We use $N(\beta\le \sigma, T)$ and $N(\beta> \sigma, T)$ to denote the count of zeros $\rho=\beta+i\gamma$ of $\Rzeta(s)$ contained in $(-\infty,\sigma]\times(0,T]$ and $(\sigma,+\infty)\times(0,T]$ respectively counting with multiplicities. The sum of the two is denoted by $N_{\Rzeta}(T)$. $N(T)$ and $N_0(T)$ denote the number of zeros of $\zeta(s)$ on the critical strip and in the critical line, respectively, as usual. We use $\Orden^*(h)$ to denote a quantity $R$ bounded by $|R|\le |h|$. 

\section{Main Results.}\label{S:main}
Applying Littlewood's lemma to $\Rzeta(s)$ in the rectangle $[\sigma,4]\times[t_0,T]$ yields 
\begin{theorem}\label{T:firstLittlewood}
Let $1\le t_0$ and  $\sigma\le1$ be fixed, then for $T\to+\infty$
\begin{equation}\label{E:firstLittlewood}
\int_{t_0}^T\log|\Rzeta(\sigma+it)|\,dt=2\pi\sum_{\substack{\beta>\sigma\\0<\gamma\le T}}(\beta-\sigma)+\Orden_\sigma(\log T).
\end{equation}
\end{theorem}
With the combed function $s\pi^{-s/2}\Gamma(s/2)\Rzeta(s)$, we can apply Littlewood's lemma on the rectangle $[\sigma_0,\sigma]\times[t_0,T]$ to get
\begin{theorem}\label{T:secondLittlewood}
Let $1\le t_0$ and $\sigma\le 10$ fixed, then  for $T\to+\infty$ we have
\begin{equation}\label{E:thm2}
2\pi\sum_{\substack{\beta\le \sigma\\0<\gamma\le T}}(\sigma-\beta)=
\sigma\Bigl(\frac{T}{2}\log\frac{T}{2\pi}-\frac{T}{2}\Bigr)+\frac{T}{2}\log 2+
\int_{t_0}^T\log|\Rzeta(\sigma+it)|\,dt+\Orden_\sigma(T^{20/21}).
\end{equation}
\end{theorem}

In some cases, we have independent information on the integral of  $\log|\Rzeta(\sigma+t)|$:
\begin{theorem}\label{T:logrzeta}
Given $\sigma_1>3$, there is a constant $C$ such that for $\sigma\ge \sigma_1$ and $1\le t_0<T$ we have 
\begin{equation}
\Bigl|\int_{t_0}^T\log|\Rzeta(\sigma+it)|\,dt\Bigr|\le C 2^{-\sigma}.
\end{equation}
\end{theorem}

\section{Applications of the main Theorems}\label{S:results}

\begin{proposition}\label{P:numberzeros}
The number of zeros $\rho=\beta+i\gamma$ of $\Rzeta(s)$ with $0<\gamma<T$ is 
\begin{equation}
N_{\Rzeta}(T)=\frac{T}{4\pi}\log\frac{T}{2\pi}-\frac{T}{4\pi}+\Orden(T^{20/21}). 
\end{equation}
\end{proposition}
\begin{proof}
By Theorems \ref{T:secondLittlewood} and \ref{T:logrzeta} we have for 
$\sigma=4$ and $\sigma=5$ for $T\to+\infty$
\begin{equation}\label{E:Partial}
2\pi\sum_{\substack{\beta\le\sigma\\0<\gamma\le T}}(\sigma-\beta)=
\sigma\Bigl(\frac{T}{2}\log\frac{T}{2\pi}-\frac{T}{2}\Bigr)+\frac{T}{2}\log 2+\Orden(T^{20/21}).\end{equation}
Except for a finite number, the zeros of $\Rzeta(s)$ with $\gamma>0$ satisfies $\beta\le 3$ as shown in \cite{A100}*{Cor.~14}.
Hence 
\[ N_{\Rzeta}(T)=\sum_{\substack{\beta\le 5\\0<\gamma\le T}}(5-\beta)-\sum_{\substack{\beta\le4\\0<\gamma\le T}}(4-\beta)+\Orden(1).\]
The result is obtained by subtracting the above results.
\end{proof}

\begin{remark}\label{R:improve}
By a different and more direct reasoning, we have proved in \cite{A185} the more precise result
\[N_{\Rzeta}(T)=\frac{T}{4\pi}\log\frac{T}{2\pi}-\frac{T}{4\pi}-\frac12\sqrt{\frac{T}{2\pi}}+\Orden(T^{2/5}\log T). \]
\end{remark}

\begin{proposition}\label{P:sumbeta}
For $T\to+\infty$ and denoting by $\rho=\beta+i\gamma$ the zeros of $\Rzeta(s)$, we have 
\begin{equation}
\sum_{0<\gamma\le T}\beta=-\frac{T}{4\pi}\log 2+\Orden(T^{20/21}).
\end{equation}
\end{proposition}
\begin{proof}
For $\sigma=4$, for example, we have 
\[\sum_{\substack{\beta\le4\\0<\gamma\le T}}(\sigma-\beta)=\sigma N_{\Rzeta}(T)-\sum_{0<\gamma\le T}\beta=\sigma\Bigl(\frac{T}{4\pi}\log\frac{T}{2\pi}-\frac{T}{4\pi}\Bigr)-\sum_{0<\gamma\le T}\beta+\Orden(T^{20/21}).\]
Comparing with \eqref{E:Partial} we get our result.
\end{proof}

\begin{lemma}\label{L:firstLemma}
Let $\sigma'<\sigma\le 10$ then 
\begin{equation}
\begin{aligned}
(\sigma-\sigma')N(\beta\le \sigma,T)&\ge(\sigma-\sigma')\Bigl(\frac{T}{4\pi}\log\frac{T}{2\pi}-\frac{T}{4\pi}\Bigr)+\\
&+\frac{1}{2\pi}\int_{t_0}^T\log|\Rzeta(\sigma+it)|\,dt-\frac{1}{2\pi}\int_{t_0}^T\log|\Rzeta(\sigma'+it)|\,dt+\orden(T).
\end{aligned}
\end{equation}
\end{lemma}
\begin{proof}
by  Theorem \ref{T:secondLittlewood} we have 
\begin{align*}
\sum_{\substack{\beta\le\sigma\\0<\gamma\le T}}(\sigma-\beta)&=
\sigma\Bigl(\frac{T}{4\pi}\log\frac{T}{2\pi}-\frac{T}{4\pi}\Bigr)+\frac{T\log2}{4\pi}+\frac{1}{2\pi}\int_{t_0}^T\log|\Rzeta(\sigma+it)|\,dt+\orden(T),\\
\sum_{\substack{\beta\le\sigma'\\0<\gamma\le T}}(\sigma'-\beta)&= 
\sigma'\Bigl(\frac{T}{4\pi}\log\frac{T}{2\pi}-\frac{T}{4\pi}\Bigr)+\frac{T\log2}{4\pi}+\frac{1}{2\pi}\int_{t_0}^T\log|\Rzeta(\sigma'+it)|\,dt+\orden(T).
\end{align*}
We subtract the two equations. Since  $\sigma'<\sigma$ we have 
\[\sum_{\substack{\beta\le\sigma\\0<\gamma\le T}}(\sigma-\beta)-
\sum_{\substack{\beta\le\sigma'\\0<\gamma\le T}}(\sigma'-\beta)=
\sum_{\substack{\beta\le\sigma\\0<\gamma\le T}}(\sigma-\beta)-
\sum_{\substack{\beta\le\sigma\\0<\gamma\le T}}(\sigma'-\beta)+
\sum_{\substack{\sigma'<\beta\le\sigma\\0<\gamma\le T}}(\sigma'-\beta)\]
\[=(\sigma-\sigma')N(\beta\le\sigma,T)+\sum_{\substack{\sigma'<\beta\le\sigma\\0<\gamma\le T}}(\sigma'-\beta)
\le (\sigma-\sigma')N(\beta\le\sigma,T).\]
Therefore,
\begin{align*}
(\sigma-\sigma')N(\beta\le \sigma,T)&\ge(\sigma-\sigma')\Bigl(\frac{T}{4\pi}\log\frac{T}{2\pi}-\frac{T}{4\pi}\Bigr)+\frac{1}{2\pi}\int_{t_0}^T\log|\Rzeta(\sigma+it)|\,dt\\&-\frac{1}{2\pi}\int_{t_0}^T\log|\Rzeta(\sigma'+it)|\,dt+\orden(T).\qedhere
\end{align*}
\end{proof}

\begin{lemma}\label{L:secondLemma}
For $\sigma'<\sigma\le 10$ we have 
\begin{equation}
(\sigma-\sigma')N(\beta\le\sigma,T)\ge(\sigma-\sigma')\Bigl(\frac{T}{4\pi}\log\frac{T}{2\pi}-\frac{T}{4\pi}\Bigr)-\frac{1}{2\pi}\int_{t_0}^T\log|\Rzeta(\sigma'+it)|\,dt +\orden(T).
\end{equation}
\end{lemma}

\begin{proof}
By Lemma \ref{L:firstLemma} and Theorem \ref{T:firstLittlewood} we have 
for any $\sigma''$ with $\sigma>\sigma''>\sigma'$
\begin{align*}
(\sigma''-\sigma')N(\beta\le\sigma,T)&\ge(\sigma''-\sigma')\Bigl(\frac{T}{4\pi}\log\frac{T}{2\pi}-\frac{T}{4\pi}\Bigr)+\sum_{\substack{\beta>\sigma''\\0<\gamma\le T}}(\beta-\sigma'')\\&-\frac{1}{2\pi}\int_{t_0}^T\log|\Rzeta(\sigma'+it)|\,dt +\orden(T).
\end{align*}
We have (by Proposition \ref{P:numberzeros})
\[\sum_{\substack{\beta>\sigma''\\0<\gamma\le T}}(\beta-\sigma'')\ge
\sum_{\substack{\beta>\sigma\\0<\gamma\le T}}(\beta-\sigma'')\ge 
(\sigma-\sigma'')N(\beta>\sigma,T)\]\[=(\sigma-\sigma'')\Bigl(\frac{T}{4\pi}\log\frac{T}{2\pi}-\frac{T}{4\pi}-N(\beta\le\sigma,T)+\orden(T)\Bigr).\]
Therefore,
\begin{align*}
(\sigma''-\sigma')N(\beta\le\sigma,T)&+(\sigma-\sigma'')N(\beta\le\sigma,T)\\&\ge(\sigma-\sigma')\Bigl(\frac{T}{4\pi}\log\frac{T}{2\pi}-\frac{T}{4\pi}\Bigr)-\frac{1}{2\pi}\int_{t_0}^T\log|\Rzeta(\sigma'+it)|\,dt +\orden(T).\qedhere
\end{align*}
\end{proof}

\begin{lemma}\label{L:Jensen}
Let $f\colon[0,+\infty)\to\C$ integrable on any compact set and such that \[\lim_{T\to+\infty}\frac{1}{T}\int_0^T|f(t)|^2\,dt=+\infty.\] 
Given  $t_0$ there is some $T_0$ such that for $T\ge T_0$ we have
\begin{equation}
\int_{t_0}^T\log|f(t)|\,dt\le \frac{T}{2}\log\Bigl(\frac{1}{T}\int_0^T|f(t)|^2\,dt\Bigr).
\end{equation}
\end{lemma}
\begin{proof}
Since $\log x$ is a concave function by  Jensen's inequality \cite{Rudin}*{p.~61} for $T>t_0$  we have
\[\frac{1}{T-t_0}\int_{t_0}^T\log|f(t)|^2\,dt\le\log\Bigl(\frac{1}{T-t_0}\int_{t_0}^T|f(t)|^2\,dt\Bigr).\]
Let us define 
\[A:=\int_{t_0}^T|f(t)|^2\,dt\le B:= \int_{0}^T|f(t)|^2\,dt.\]
The function 
$x\mapsto \frac{x}{2}\log \frac{B}{x}$
is increasing for $x\in[T-t_0,T]$, since the derivative is 
$\frac12\log\frac{B}{x}-\frac12$, which is
positive for $B>e x$. So, we only need that $B>eT$ and this happens for $T\ge T_0$ by hypothesis. Then 
\[\int_{t_0}^T\log|f(t)|\,dt\le\frac{T-t_0}{2}\log\frac{A}{T-t_0} \le\frac{T-t_0}{2}\log\frac{B}{T-t_0} \le\frac{T}{2}\log\frac{B}{T}.\]
This is our inequality.
\end{proof}

\begin{proposition}\label{P:Npar}
There is some $T_0$ such that for $T>T_0$ we have
\begin{equation}
N(\beta\le\tfrac12,T)\ge \frac{35T}{396\pi}+\frac{5}{11\pi}\int_{t_0}^T\log|\Rzeta(\tfrac12+it)|\,dt.
\end{equation}
\end{proposition}

\begin{proof}
Applying  Lemma \ref{L:firstLemma} with $\sigma=1/2$ and $\sigma'<1/2$ yields 
\begin{align*}
(\tfrac12-\sigma')N(\beta\le\tfrac12,T)&\ge(\tfrac12-\sigma')\Bigl(\frac{T}{4\pi}\log\frac{T}{2\pi}-\frac{T}{4\pi}\Bigr)+\\
&+\frac{1}{2\pi}\int_{t_0}^T\log|\Rzeta(\tfrac12+it)|\,dt-\frac{1}{2\pi}\int_{t_0}^T\log|\Rzeta(\sigma'+it)|\,dt+\orden(T).
\end{align*}
In \cite{A101} we proved that for $\sigma'\le 1/4$ we have 
\[\frac{1}{T}\int_0^T|\Rzeta(\sigma'+it)|^2\,dt=\frac{2}{(1-2\sigma')(3-2\sigma')}\Bigl(\frac{T}{2\pi}\Bigr)^{\frac12-\sigma'}+\Orden(T^{\frac14-\sigma'}).\]
This and Lemma \ref{L:Jensen} yields for $\sigma\le1/4$ and $T$ big enough
\[\int_{t_0}^T\log|\Rzeta(\sigma'+it)|\,dt\le(\tfrac12-\sigma')\frac{T}{2}\log\frac{T}{2\pi}-\frac{T}{2}\log\frac{(1-2\sigma')(3-2\sigma')}{2}+\Orden(T^{3/4}).\]
Hence,
\begin{align*}
(\tfrac12-\sigma')N(\beta\le\tfrac12,T)&\ge-(\tfrac12-\sigma')\frac{T}{4\pi}+\\
&+\frac{1}{2\pi}\int_{t_0}^T\log|\Rzeta(\tfrac12+it)|\,dt+\frac{T}{4\pi}\log\frac{(1-2\sigma')(3-2\sigma')}{2}+\orden(T).
\end{align*}
The function 
\[f(\sigma')=-1+\frac{1}{\frac12-\sigma'}\log\frac{(1-2\sigma')(3-2\sigma')}{2},\]
has a maximum near $\sigma'=-3/5$ where $f(-3/5)>7/18$. Therefore, taking $\sigma'=-3/5$ we get
\[\frac{11}{10}N(\beta\le\tfrac12,T)\ge \frac{7T}{72\pi}+\frac{1}{2\pi}\int_{t_0}^T\log|\Rzeta(\tfrac12+it)|\,dt+\orden(T).\]
Since $f(-3/5)>7/18$ the term $\orden(T)$ can be eliminated for $T\ge T_0$, taking $T_0$ big enough. 
\end{proof}

\begin{corollary}\label{C:leftzeros}
There exists some $T_0$ such that for $T\ge T_0$ we have 
\[N(\beta\le\tfrac12,T)\ge \frac{35T}{396\pi}+\frac{5}{11}\sum_{\substack{\beta>1/2\\0<\gamma\le T}}(\beta-\tfrac12).\]
\end{corollary}

\begin{proof}
In Proposition \ref{P:Npar} substitute the integral by its value in Theorem \ref{T:firstLittlewood}. The error term $\Orden(\log T)$ is eliminated due to the strict inequality $f(-3/5)>7/18$.
\end{proof}
\section{Connection with the zeros of \texorpdfstring{$\zeta(s)$}{zeta}}

We have $\zeta(\frac12+it)=e^{-i\vartheta(t)}Z(t)$ where $\vartheta(t)$ and $Z(t)$ are real analytic functions. These are connected to $\Rzeta(s)$ by 
\[Z(t)=2\Re\bigl\{e^{i\vartheta(t)}\Rzeta(\tfrac12+it)\bigr\},\]
(see \cite{A166}). In \cite{A66} we show that $\Rzeta(\frac12+it)=e^{-i\omega(t)}g(t)$ where $\omega(t)$ and $g(t)$ are real analytic functions. It follows that a point $\frac12+it$ with $t\in\R$ is a zero of $\zeta(s)$ if and only if it is a zero of $\Rzeta(s)$ or if 
\[\cos\bigl(\vartheta(t)-\omega(t)\bigr)=0.\]
Also, it is shown in \cite{A66} that for $T>0$ we have 
\[\omega(T)=2\pi N(\beta\ge1/2,T)+\Orden(\log T).\]
If there were no zeros to the right of the critical line, $N(\beta\ge1/2,T)=0$, then the function $\vartheta(t)-\omega(t)$ will increase for $t\in[0,T]$ from $0$ to $\approx\vartheta(T)$ and the number of zeros of $\zeta(s)$ in the critical line from $t=0$ to $t=T$  will be $\vartheta(T)/\pi$.  Almost all zeros of the zeta function will be on the critical line. But this does not appear to be true, our (limited) computation of zeros \cite{A172} points to the approximate equality $N(\beta\ge1/2,T)\approx \vartheta(T)/6\pi$. 

\section{Proof of Theorem 3}\label{S:Th3}

\begin{lemma}\label{L:zetaRzeta}
For $s=\sigma+it$ with $\sigma\ge2$, $t\ge0$ and $|s|>2\pi e^2$ we have
\begin{equation}
\Rzeta(s)=\zeta(s)(1+W(s)),\qquad |W(s)|\le C\Bigl(\frac{|s|}{2\pi e}\Bigr)^{\frac{1-\sigma}{2}},
\end{equation}
where $C$ is an absolute constant.
\end{lemma}
\begin{proof}
In \cite{A100}*{Thm.~12} it is proved that for $s=\sigma+it$ with $t>0$,  $\sigma>0$ and such that $|s|>2\pi e^2$ we have
\[\Rzeta(s)=\sum_{n=1}^\ell\frac{1}{n^s}+R(s),\qquad |R(s)|\le C |s/2\pi e|^{-\sigma/2},\]
where, $\xi=\sqrt{\frac{s}{2\pi i}}=\xi_1+i\xi_2$ with $0<-\xi_2<\xi_1$ and 
$\ell=\lfloor \xi_1-\xi_2\rfloor$, and $C$ is an absolute constant.

Therefore, under these conditions, but assuming $\sigma>1$
\[\Rzeta(s)=\zeta(s)-\sum_{n=\ell+1}^\infty\frac{1}{n^s}+R(s).\]
Notice that $|\xi_1|=\xi_1$ and  $|\xi_2|=-\xi_2$, so that
\[|\xi|-1\le|\xi_1|+|\xi_2|-1< \ell=\lfloor \xi_1-\xi_2\rfloor\le|\xi_1|+|\xi_2|\le \sqrt{2}|\xi|=\sqrt{|s|/\pi}.\]
Therefore,
\[\Bigl|\sum_{n=\ell+1}^\infty\frac{1}{n^s}\Bigr|\le \sum_{n>|\xi|}\frac{1}{n^\sigma}\le |\xi|^{-\sigma}+\int_{|\xi|}^\infty x^{-\sigma}\,dx=
(|s|/2\pi)^{-\sigma/2}+\frac{1}{\sigma-1}(|s|/2\pi)^{(1-\sigma)/2}.\]
It follows that for $\sigma\ge2$ 
\[\Bigl|R(s)-\sum_{n=\ell+1}^\infty\frac{1}{n^s}\Bigr|\le C |s/2\pi e|^{-\sigma/2}+
(|s|/2\pi)^{-\sigma/2}+\frac{1}{\sigma-1}(|s|/2\pi)^{(1-\sigma)/2}\le C|s/2\pi e|^{(1-\sigma)/2}.\]
Therefore,
\[\Rzeta(s)=\zeta(s)+V(s)=\zeta(s)(1+W(s)),\qquad |W(s)|\le C\zeta(2)|s/2\pi e|^{(1-\sigma)/2}.\qedhere\]
\end{proof}

\begin{proposition}\label{P:exactR}
There exist $\sigma_0>3$, $t_0>1$ and a constant $C$ such that for $T>t_0$ and $\sigma>\sigma_0$, we have
\begin{equation}
\Bigl|\int_{t_0}^T\log|\Rzeta(\sigma+it)|\,dt\Bigr|\le C2^{-\sigma},
\end{equation}
\end{proposition}

\begin{proof}
We have
\[\int_{t_0}^T\log|\Rzeta(\sigma+it)|\,dt=\int_{t_0}^T\log|\zeta(\sigma+it)|\,dt+\int_{t_0}^T|\log(1+W(\sigma+it))|\,dt.\]
The integral in the zeta function is treated easily
\begin{align*}
\int_{t_0}^T\log|\zeta(\sigma+it)|\,dt&=\Re\int_{t_0}^T\sum_{n=2}^\infty\frac{\Lambda(n)}{\log n}\frac{1}{n^{s}}\,dt=\int_{t_0}^T\sum_{n=2}^\infty\frac{\Lambda(n)}{\log n}\frac{\cos(t\log n)}{n^{\sigma}}\,dt\\
&=\sum_{n=2}^\infty\frac{\Lambda(n)}{\log n}\frac{\sin(T\log n)-\sin(t_0\log n)}{n^{\sigma}\log n}.
\end{align*}
So, in  absolute value 
\[\Bigl|\int_{t_0}^T\log|\zeta(\sigma+it)|\,dt\Bigr|\le 2^{-\sigma}\sum_{n=2}^\infty \frac{2\Lambda(n)}{(n/2)^\sigma\log^2n}\le 2^{-\sigma}\sum_{n=2}^\infty \frac{2\Lambda(n)}{(n/2)^3\log^2n}.\]

The function $1+W(s)$ defined by $\Rzeta(s)=\zeta(s)(1+W(s))$ is analytic and is not equal to $0$ on the line $\sigma+it$ with $\sigma>3$ fixed and $t\ge t_0$. 
Since $|W(s)|\le C(|s|/2\pi e)^{(1-\sigma)/2}$, there is some $t_1$ such that 
for $t>t_1$ we have $|W(s)|\le 1/4$.  The function $\log(1+W(s))$ is well defined,  with the imaginary part $\le \pi/2$ in absolute value for $t>t_1$. 
We also have for  $t\ge t_1$ that 
\[|\log(1+W(s))|\le 2C(|s|/2\pi e)^{(1-\sigma)/2}.\]
Then this inequality is also true for any $t\ge t_0$, substituting, if needed, 
$C$ by a large constant. 

Hence, 
\[\Bigl|\int_{t_0}^T\log(1+W(t))\,dt\Bigr|\le C'\int_{t_0}^T|s/2\pi e|^{\frac{1-\sigma}{2}}\,dt=C'(2\pi e)^{\frac{\sigma-1}{2}}\int_1^T(\sigma^2+t^2)^{\frac{1-\sigma}{4}}\,dt\]
\[\le C'\sigma^{1+\frac{1-\sigma}{2}}(2\pi e)^{\frac{\sigma-1}{2}}\int_{1/\sigma}^{T/\sigma}(1+u^2)^{\frac{1-\sigma}{4}}\,du\le C'\sigma\Bigl(\frac{2\pi e}{\sigma}\Bigr)^{\frac{\sigma-1}{2}}\int_0^{\infty}(1+u^2)^{\frac{1-\sigma}{4}}\,du. \]
Elementary computations and Euler-MacLaurin approximation to the Gamma function yields for $\sigma>3$
\[=C'\sigma\Bigl(\frac{2\pi e}{\sigma}\Bigr)^{\frac{\sigma-1}{2}}\frac{\sqrt{\pi}\Gamma(\frac{\sigma-3}{4})}{2\Gamma(\frac{\sigma-1}{4})}\le 
C''\frac{\sigma^{3/2}}{\sigma-3}\Bigl(\frac{2\pi e}{\sigma}\Bigr)^{\frac{\sigma-1}{2}}.\]
We end noticing that for $\sigma>\sigma_0>3$, there is a constant $C$ only depending on $\sigma_0$ such that
\[2^{-\sigma}+\frac{\sigma^{3/2}}{\sigma-3}\Bigl(\frac{2\pi e}{\sigma}\Bigr)^{\frac{\sigma-1}{2}}\le C 2^{-\sigma}.\qedhere\]
\end{proof}

\section{Proof of Theorem 1}

We will use Littlewood and Backlund's lemmas, which we state as reference.

\begin{lemma}[Littlewood]
Let $f\colon\Omega\to\C$ be a holomorphic function, $R=[a,b]\times[c,d]\subset\Omega$ a closed rectangle contained in $\Omega$, let us denote by $A=a+ic$, $B=b+ic$, $C=b+id$, and 
$D=a+id$ its vertices. Assume that $f$ do not vanish on the sides 
$AB$, $BC$ and $CD$. Define $\arg f(z)$ continuously along this three sides then 
\begin{equation}
\begin{aligned}
2\pi\sum_{\substack{\beta\ge a\\c<\gamma\le d}}(\beta-a)&=\int_c^d\log|f(a+iy)|\,dy-\int_c^d \log|f(b+iy)|\,dy\\ &-\int_a^b\arg f(x+ic)\,dx+\int_a^b\arg f(x+id)\,dx.
\end{aligned}
\end{equation}
where $\rho=\beta+i\gamma$ run through the zeros of $f(z)$ in the rectangle counted with multiplicities. 
\end{lemma}

\begin{lemma}[Backlund]\label{L:B}
Let  $f$ be holomorphic in the disc $|z-a|\le R$. 
Let  $|f(z)|\le M$ for  $|z-a|\le R$. Let  $b$ be a point in the interior of the disc
$0<|b-a|<R$. Assume that  $f$ does not vanish on the segment $[a,b]$, then 
\begin{equation}\label{E:backlund}
\Bigl|\Re\frac{1}{2\pi i}\int_a^b\frac{f'(z)}{f(z)}\,dz\Bigr|\le
\frac{1}{2}\Bigl(\log \frac{M}{|f(a)|}\Bigr)\Bigl(\log\frac{R}{|b-a|}\Bigr)^{-1}.
\end{equation}
\end{lemma}
\begin{remark}
Littlewood's lemma can be found in many sources, for example \cite{T}*{section~9.9} or the original \cite{L}. Backlund's lemma is applied without explicit mention in  Backlund \cite{B}. When stated, an extra term $1/2$ is frequently added on the right side. It is proved by means of the Jensen formula \cite{Rudin}*{section 15.16}. An independent complete proof is given in \cite{A185}.
\end{remark}

\begin{proof}[Proof of Theorem \ref{T:firstLittlewood}]\label{Proof 1Littl}
The truth of  \eqref{E:firstLittlewood} does not depend on the particular value of $t_0$. Changing $t_0$ to some other value changes the left-hand side of \eqref{E:firstLittlewood} into a constant that is absorbed into the error term.
So we can assume that $\Rzeta(s)\ne0$ for $\Im(s)=t_0$, and that $t_0>32\pi$. 

We can also assume that $\Rzeta(s)\ne0$ for  $\Im(s)=T$. In the other case, we  may prove the proposition for $T_n>T$ with $\lim_n T_n=T$, and obtain the proposition in the general case by taking the limits of the resulting equation for $n\to\infty$.  

For $\Re(s)=4$ and $t\ge 32\pi$ we have $|\Rzeta(s)-1|<\frac12$  by  Proposition 6 in \cite{A173}. Therefore,  $\Rzeta(4+it)$ do not vanish for $t\ge32\pi$, and we are in a condition to apply Littlewood's lemma to $\Rzeta(s)$ in the rectangle $[\sigma, 4]\times[t_0,T]$.

This yields
\begin{align*}\text{\checked}\quad
2\pi\sum_{\substack{\beta>\sigma\\t_0<\gamma\le T}}(\beta-\sigma)=
\int_{t_0}^T\log|\Rzeta(\sigma+it)|\,dt-\int_{t_0}^T\log|\Rzeta(4+it)|\,dt\\-
\int_{\sigma}^4\arg\Rzeta(x+it_0)\,dx+\int_{\sigma}^4\arg\Rzeta(x+iT)\,dx,
\end{align*}
where $\arg\Rzeta(x+it_0)$ and $\arg\Rzeta(x+iT)$ are continuous extensions of  $\arg\Rzeta(4+it)$. Since $|\Rzeta(4+it)-1|<1/2$,  we may take $\arg\Rzeta(4+it)$
continuous and with absolute value $<\pi/2$.

On the upper side $\Im(s)=T$ we have 
\[\arg\Rzeta(x+iT)-\arg\Rzeta(4+iT)=\Re\frac{1}{ i}\int_{4+iT}^{x+iT}\frac{\Rzeta'(z)}{\Rzeta(z)}\,dz.\]
The absolute value of this integral can be bounded as usual with Backlund's lemma \ref{L:B}.  Let $D$ be the disc with center at $4+iT$ and radius 
$2(4-\sigma)$. The maximum of $|\Rzeta(s)|$ on this disc is determined by Propositions 12 and 13 in \cite{A92}, it is $\le C_\sigma T^{r}$ where $r=\max(1/2,3/2-\sigma)$. Since we assume $\sigma\le1$, we have $r\le 3/2-\sigma$. [In this argument $\sigma$ is fixed and we assume $T>1-\sigma$, for example].  Note also that $|\Rzeta(4+iT)|>1/2$. 
It follows that 
$|\arg\Rzeta(x+iT)|\le C'_\sigma\log T$ so that 
\[\Bigl|\int_{\sigma}^4\arg\Rzeta(x+iT)\,dx\Bigr|\le C''_\sigma \log T.\]
The integral of $\arg\Rzeta(x+it_0)$ is a constant, depending on $\sigma$. 

By proposition \ref{P:exactR} the integral of $\Rzeta(4+it)$ is bounded by an absolute constant. Joining all this, Littlewood's lemma yields \eqref{E:firstLittlewood}
\end{proof}

\section{Proof of Theorem \ref{T:secondLittlewood}}\label{S:Th2}

\begin{figure}[H]
\begin{center}
\includegraphics[width=0.7\hsize]{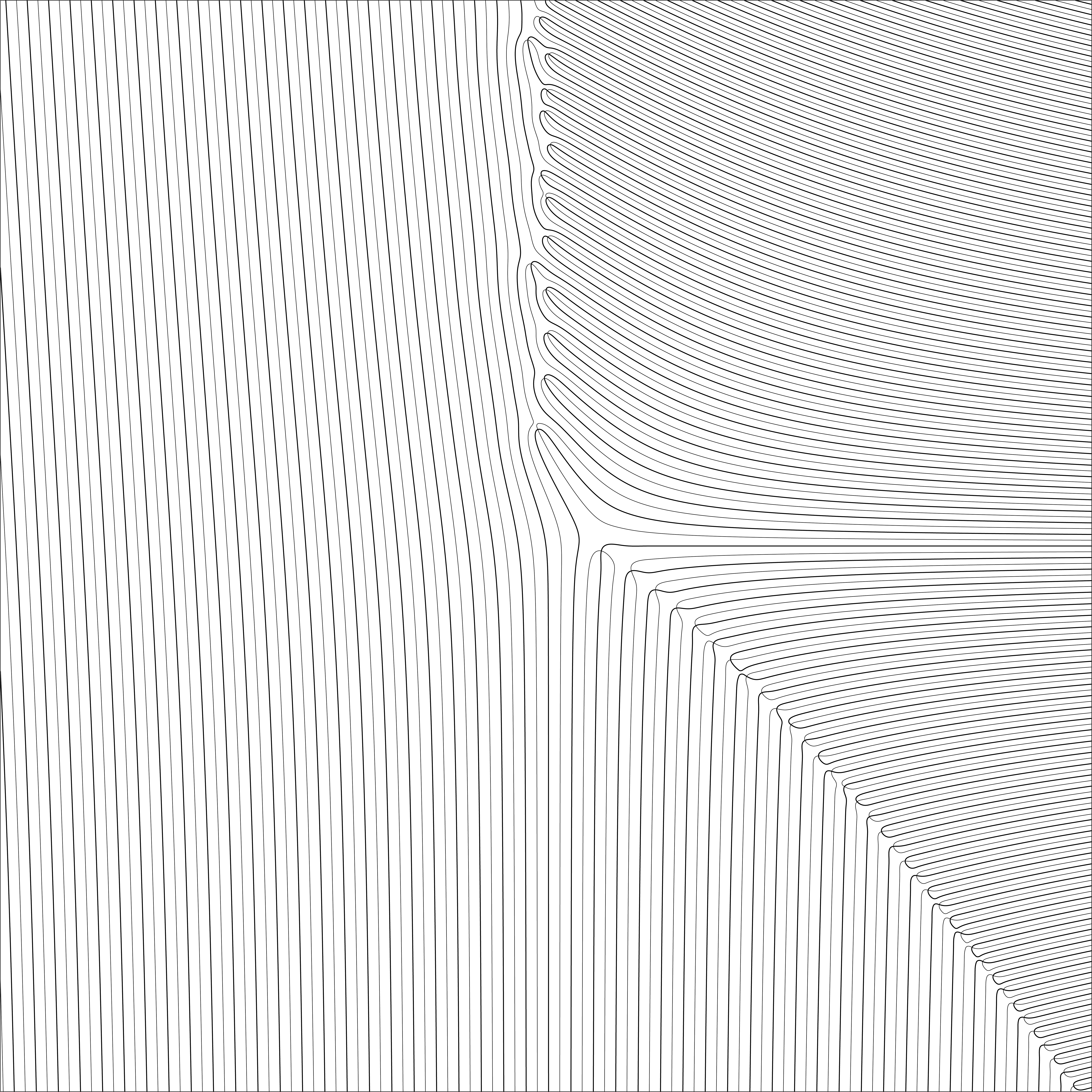}
\caption{x-ray of $F(s)=s\pi^{-s/2}\Gamma(s/2)\Rzeta(s)$ on $(-105,105)^2$}
\end{center}
\end{figure}

To prove Theorem \ref{T:secondLittlewood} Siegel \cite{Siegel} consider a function 
\[g(s)=\pi^{-\frac{s+1}{2}}e^{-\frac{\pi i s}{4}}\Gamma(\tfrac{s+1}{2})\Rzeta(s).\]
The election of the factor of $\Rzeta(s)$ is to make the behavior of $g(\sigma_0+it)$  the simplest possible, for our election of $\sigma_0$. Here $\sigma_0$ is selected so that all zeros to height $T$ of $\Rzeta(s)$ satisfies $\Re(\rho)>\sigma_0$.  Nevertheless, this election makes $g$ meromorphic with poles at $-(2n+1)$. This is a problem to apply Littlewood's lemma. Therefore, we prefer to instead consider the entire function 
\begin{equation}
F(s)=s\pi^{-s/2}\Gamma(\tfrac{s}{2})\Rzeta(s).
\end{equation}

\begin{remark}
The x-ray of Siegel's function $g(s)$ is very nice.
\begin{figure}[H]
\begin{center}
\includegraphics[width=0.7\hsize]{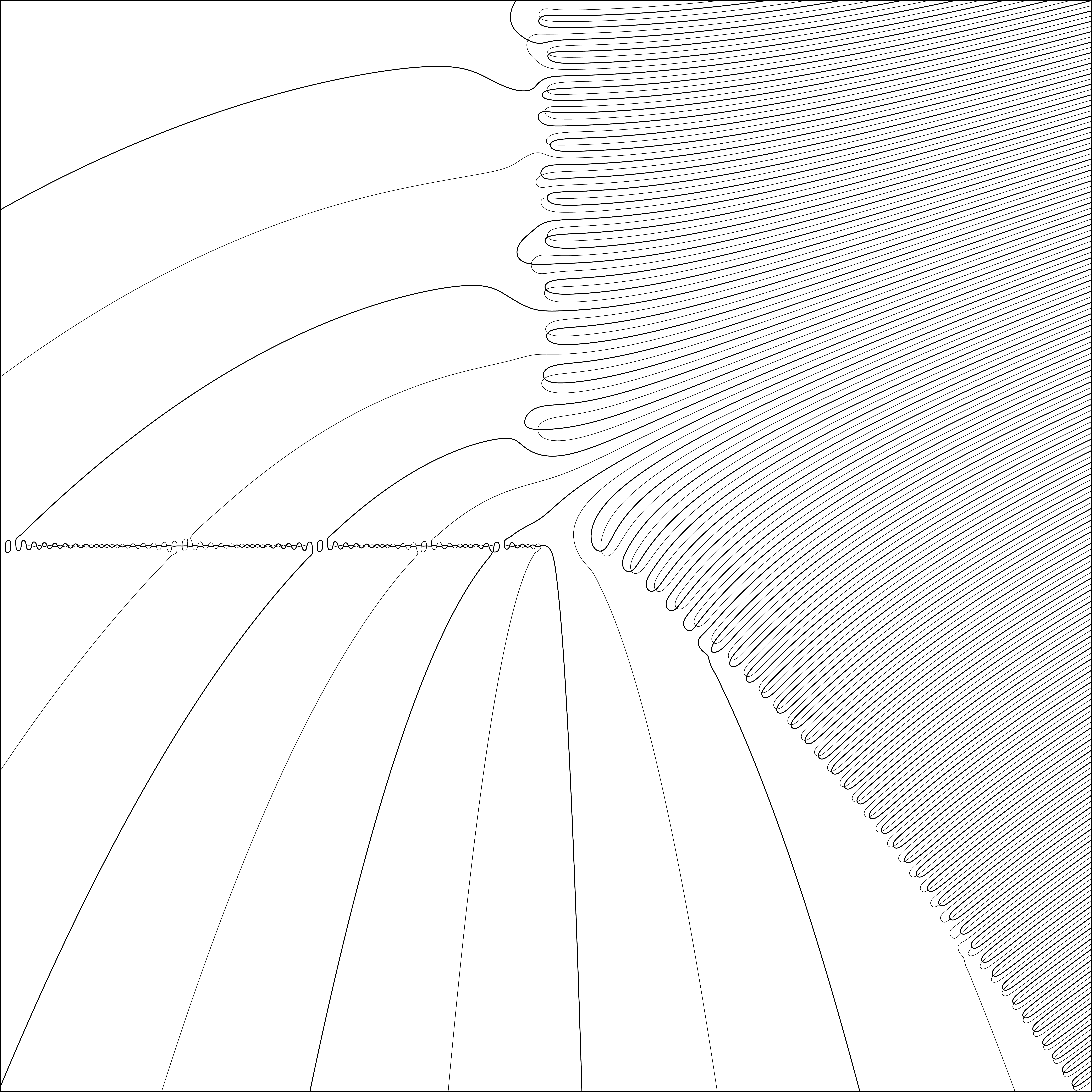}
\caption{x-ray of $g(s)=\pi^{-\frac{s+1}{2}}e^{-\frac{\pi i s}{4}}\Gamma(\tfrac{s+1}{2})\Rzeta(s)$ on $(-105,105)^2$}
\label{default1}
Since it has poles at $-1$, $-3$, \dots\  and zeros at $-2$, $-4$, \dots the x-rays contains cycles. We can see their delicate game in the attached detail figure. 
\begin{figure}[H]
\begin{center}
\includegraphics[width=\hsize]{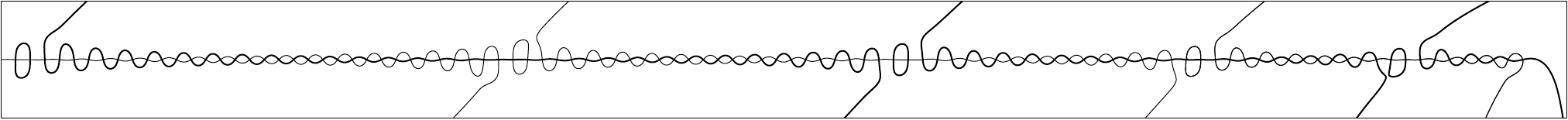}
\label{default2}
\end{center}
\end{figure}
\begin{figure}[H]
\begin{center}
\includegraphics[width=0.5\hsize]{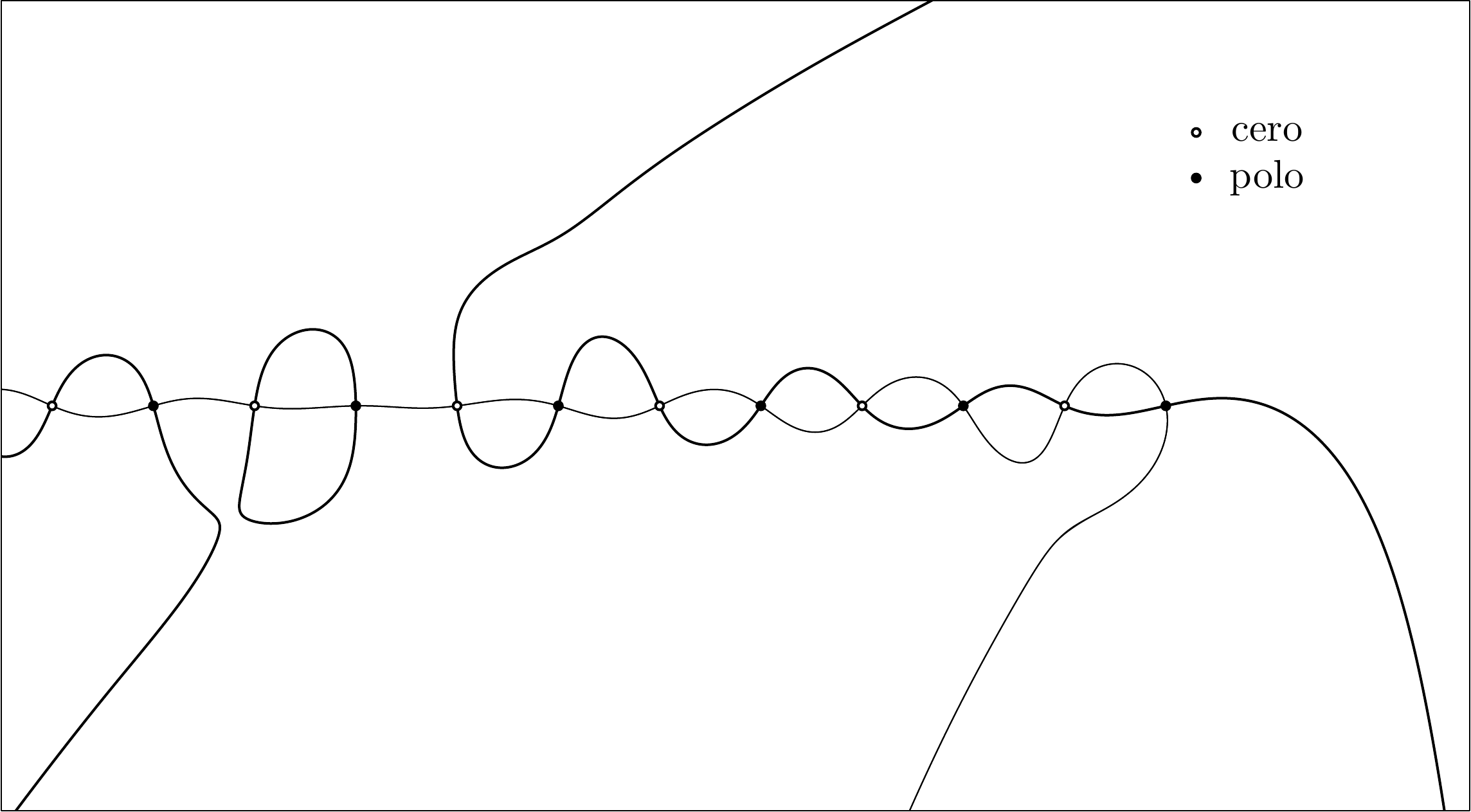}
\caption{xray of $g(s)$ details}
\label{default3}
\end{center}
\end{figure}

\end{center}
\end{figure}

\end{remark}

\textbf{Notation:} In this section,  we will consider the values of $\Rzeta(s)$ at the points $s=\sigma+it$ with $\sigma\le 10$ and $t\ge0$. For such values of $s$, we consider $\eta=\sqrt{\frac{s-1}{2\pi i}}$, selecting the root  $\eta=\eta_1+i\eta_2$ with $\eta_1+\eta_2\ge0$. When $\sigma\le 1$ this means $0\le \eta_2\le \eta_1$. In this case 
$|e^{2\pi i \eta}|\le 1$, with strict inequality  for  $\sigma<1$. 

We want a value of $\sigma$ such that all zeros of $\Rzeta(s)$ to height $T$ satisfies $\sigma\le\Re(\rho)$. This is obtained by means of the next Theorem proved in  \cite{A193}*{Th.~9}.
\begin{theorem}\label{T:RzetaAprox}
There exist constants $A$ and $t_0>1$ such that for $s$ in the closed set
\[\Omega=\{s\in\C\colon t\ge t_0,\quad 1-\sigma\ge t^{3/7}\},\]
we have 
\begin{equation}\label{E:withUT3}
\Rzeta(s)=-\chi(s)\eta^{s-1}e^{-\pi i \eta^2}\frac{\sqrt{2}e^{3\pi i/8}\sin\pi\eta}{2\cos2\pi\eta}(1+U(s)),\qquad |U(s)|\le At^{-\frac{1}{21}},
\end{equation}
where $\eta=\sqrt{(s-1)/2\pi i}$ satisfies $\Re(\eta)+\Im(\eta)>0$. 
\end{theorem}

\begin{lemma}\label{L:Jexpres}
Let $t_0$ be a constant greater than the constant appearing in Theorem \ref{T:RzetaAprox}. 
For $s=\sigma+it$ with $t\ge t_0$ and $1-\sigma\ge t^{3/7}$, we have
\[F(s):=s\pi^{-s/2}\Gamma(s/2)\Rzeta(s)=J(s)(1+U(s)),\]
where 
\begin{equation}\label{E:logJ}
\checked\log J(s)=\frac{\pi i s}{4}+\log s-\frac{1}{2}\log\frac{1-s}{2\pi }-\frac{3\pi i}{8}+\pi i\eta+\log\frac{1-e^{2\pi i \eta}}{1+e^{4\pi i \eta}}+\Orden(|s|^{-1}).\end{equation}
\end{lemma}
\begin{proof}
Since we assume that $s$ satisfies the conditions in Theorem \ref{T:RzetaAprox} we have 
\eqref{E:withUT3}. Therefore, by definition,  $F(s)=J(s)(1+U(s))$, where 
\begin{align*}
J(s)&=-s\pi^{-s/2}\Gamma(s/2)\chi(s)\eta^{s-1}e^{-\pi i \eta^2}\frac{\sqrt{2}e^{3\pi i/8}\sin\pi\eta}{2\cos2\pi\eta}\\
&=-s\pi^{(s-1)/2}\Gamma(\tfrac{1-s}{2})\eta^{s-1}e^{-\pi i \eta^2}\frac{\sqrt{2}e^{3\pi i/8}\sin\pi\eta}{2\cos2\pi\eta}\\
&=2^{-1/2}e^{-\pi i/8}s\pi^{(s-1)/2}\Gamma(\tfrac{1-s}{2})\eta^{s-1}e^{-\pi i \eta^2+\pi i \eta}\frac{1-e^{2\pi i \eta}}{1+e^{4\pi i \eta}}.
\end{align*}
We may define a continuous logarithm of $J(s)$ for these values of $s$ by noticing that $|e^{2\pi i\eta}|<1$ and therefore we may use 
\[\log\frac{1-e^{2\pi i \eta}}{1+e^{4\pi i \eta}}=-\sum_{n=1}^\infty\frac{z^n}{n}+\sum_{n=1}^\infty(-1)^n\frac{z^{2n}}{n},\qquad z=e^{2\pi i \eta}.\]
That is, it is equal to $\log(1-e^{2\pi i \eta})-\log(1+e^{4\pi i \eta})$, using in both cases the main branch of the logarithm. It follows that the imaginary part of $\log\frac{1-e^{2\pi i \eta}}{1+e^{4\pi i \eta}}$ is bounded by $\pi$. 
By the usual Euler-MacLaurin expansion, which is applicable since $\sigma<0$, 
\begin{multline*}\log J(s)=-\frac12\log2-\frac{\pi i}{8}+\log s+\frac{s-1}{2}\log \pi-\frac{s}{2}\log\frac{1-s}{2}-\frac{1-s}{2}+\frac12\log2\pi\\+
\frac{s-1}{2}\log\frac{s-1}{2\pi i}-\pi i\frac{s-1}{2\pi i}+\pi i\eta+\log\frac{1-e^{2\pi i \eta}}{1+e^{4\pi i \eta}}+\Orden(|s|^{-1}).\end{multline*}
Simplifying this yields 
\[\log J(s)=\frac{\pi i s}{4}+\log s-\frac12\log(1-s)+\frac12\log2\pi-\frac{3\pi i}{8}+\pi i\eta+\log\frac{1-e^{2\pi i \eta}}{1+e^{4\pi i \eta}}+\Orden(|s|^{-1}).\qedhere\]
\end{proof}

The next propositions will be proved later in Section \ref{S:comp}.  

\begin{proposition}\label{P:logF}
Let $t_0>1$ and $0<a<\frac12$  be  fixed real numbers. There is a function $f\colon[0,+\infty)\to[0,+\infty)$ such that $f(T)=\orden(T^{1/2})$ and such that 
for $T>t_0$,  $\sigma\le 10$ and $1-\sigma\le T^a$ we have
\begin{equation}
\begin{aligned}
\int_{t_0}^T \log|F(\sigma+it)|\,dt&=
-\frac{\pi T^2}{8}+(\sigma+1)\Bigl(\frac{T}{2}\log \frac{T}{2\pi}-\frac{T}{2}\Bigr)+\frac{T}{2}(\log 2+2\log(2\pi))\\&-\frac{\pi\sigma^2}{8}+\int_{t_0}^T\log|\Rzeta(\sigma+it)|\,dt+\Orden(f(T)).
\end{aligned}
\end{equation}
\end{proposition}

\begin{proposition}\label{P:intvert}
Let $1\le t_0<T$ with $t_0$ bigger than the constant appearing in Theorem \ref{T:RzetaAprox} and let $1-\sigma_0=T^a$ with $a=3/7$.  Then we have 
\begin{equation}
\int_{t_0}^T\log|F(\sigma_0+it)|\,dt=-\frac{\pi T^2}{8}+\frac{T}{2}\log T-\frac{T}{2}+\frac{T}{2}\log2\pi+\Orden(T^{20/21}).
\end{equation}
\end{proposition}

\begin{lemma} \label{L:underside}
Let $t_0$ and $\sigma_0$ be fixed real numbers satisfying $t_0>1$ and $\sigma_0<-1$. Assume that $F(s)$ does not vanish for $\Im(s)=t_0$ and denote by $\arg F(\sigma +it_0)$ a continuous determination of the argument, defined for $\sigma_0\le\sigma\le\sigma_1\le 10$, and such that  $|\arg F(\sigma_0+it_0)|\le C|\sigma_0|$. Then there is a constant $C_0$ (depending on $t_0$, and $C$ such that 
\[\Bigl|\int_{\sigma_0}^{\sigma_1}\arg F(\sigma+it_0)\,d\sigma\Bigr|\le C_0
|\sigma_0|^2\log|\sigma_0|.\]
\end{lemma}

\begin{proof}
We have 
\[|\arg F(\sigma+it_0)-\arg F(\sigma_0+it_0)|=\Bigl|\Re\frac{1}{i}\int_{\sigma_0+it_0}^{\sigma+it_0}\frac{F'(z)}{F(z)}\,dz\Bigr|,\]
we may bound the integral by Backlund's lemma. We take a disc $D$ with center at $10+it_0$ and radius $2(10-\sigma_0)$. To bound $F$ on $D$, notice that $D$ is contained in a disc of center $0$ and radius $C|\sigma_0|$ (where the constant $C$ here depends only on $t_0$). In \cite{A66}*{eq.(3)} we have proved that $|F(s)|\le A e^{B|s|\log |s|}$.
So that $|F(s)|\le A e^{c(t_0) |\sigma_0|\log|\sigma_0|}$. Also, we have $|F(10+it_0)|>0$ is a constant depending only on $t_0$. Therefore,
\[|\arg F(\sigma+it_0)-\arg F(\sigma_0+it_0)|\le C|\sigma_0|\log|\sigma_0|.\]
It follows that $|F(\sigma+it_0)|\le C|\sigma_0|\log|\sigma_0|$. Then 
\[\Bigl|\int_{\sigma_0}^{\sigma_1}\arg F(\sigma+it_0)\,d\sigma\Bigr|\le C
|\sigma_0|^2\log|\sigma_0|.\qedhere\]
\end{proof}

\begin{lemma}\label{L:upperside}
Let $0<a<\frac12$ a fixed real number and for $T>1$ let $1-\sigma_0=T^a$ and $\sigma_0<\sigma_1\le 10$. Then we have
\[\Im\int_{\sigma_0+iT}^{\sigma_1+iT}\log(s\pi^{-s/2}\Gamma(s/2))\,ds=\Orden(T^{2a})\]
when we take a continuous branch of the logarithm where $\Im\log(s\pi^{-s/2}\Gamma(s/2))$ at the point $s=s_0+iT$ is $\Orden(T^a)$. 
\end{lemma}

\begin{proof}
By Euler-MacLaurin expansion we have
\[\checked \Im\log\Gamma((\sigma+iT)/2)=\frac{T}{2}\log\frac{T}{2}-\frac{T}{2}+\frac{\pi}{4}(\sigma-1)+\frac{2\sigma-\sigma^2}{4T}+\Orden(T^{-1}).\]
For $\sigma_0\le\sigma\le \sigma_1\le10$ we have $\sigma^2/T^2=\Orden(T^{-1})$ so that 
\[\Im\log(\sigma+iT)=\frac{\pi}{2}-\arctan\frac{\sigma}{T}=\frac{\pi}{2}-\frac{\sigma}{T}+\Orden(T^{-1})\]
After some simplifications, we obtain for $s=\sigma+iT$
\[\text{\checked}\Im\log(s\pi^{-s/2}\Gamma(s/2))=\frac{T}{2}\log\frac{T}{2\pi}-\frac{T}{2}-\frac{\sigma^2+2\sigma}{4T}+\frac{\pi}{4}+\frac{\pi\sigma}{4}+\Orden(T^{-1}).\]
After subtracting the multiple of $2\pi$ nearest to $\frac{T}{2}\log\frac{T}{2\pi}-\frac{T}{2}$, we get a branch of the logarithm equal to 
\[\Im\log(s\pi^{-s/2}\Gamma(s/2))=-\frac{\sigma^2+\sigma}{2T}+\frac{\pi}{4}+\frac{\pi\sigma}{4}+\Orden(1)=\frac{\pi\sigma}{4}+\Orden(1).\]
Therefore, 
\[\Bigl|\Im\int_{\sigma_0+iT}^{\sigma_1+iT}\log(s\pi^{-s/2}\Gamma(s/2))\,ds\Bigr|\le \frac{\pi}{8}(\sigma_1^2-\sigma_0^2)+ C(\sigma_1-\sigma_0)\le C T^{2a}.\qedhere\]
\end{proof}

\begin{proof}[Proof of Theorem \ref{T:secondLittlewood}]
For a given $T>t_0$ define $\sigma_0$ such that $1-\sigma_0=T^a$ with $a=3/7$. For $T\ge T_0$ we will have $\sigma_0<\sigma$. We apply Littlewood's lemma in the rectangle $R=[\sigma_0,\sigma]\times[t_0,T]$ to the function $F(s)=s\pi^{-s/2}\Gamma(s/2)\Rzeta(s)$. As we have seen  in the proof  of Theorem \ref{T:firstLittlewood} we may assume that there is no zero of $F(s)$ on the lines $\Re(s)=t_0$ and $\Re(s)=T$. We take $t_0$ large enough so that Theorem \ref{T:RzetaAprox} applies with $At_0^{-1/21}<1$. This implies that $F(s)$ do not vanish on the left-hand side of the rectangle $R$. Note that the zeros of $F(s)$ and $\Rzeta(s)$ on $t>0$ are the same with the same multiplicities. 

Here, it is in line $\Re(s)=\sigma$ where we do not know the argument of $F(s)$. Hence we have
\begin{equation}\label{E:Litequation}
\begin{aligned}\text{\checked}\quad
2\pi\sum_{\substack{\beta<\sigma\\t_0<\gamma\le T}}(\sigma-\beta)=
\int_{t_0}^T\log|F(\sigma+it)|\,dt-\int_{t_0}^T\log|F(\sigma_0+it)|\,dt\\
+\int_{\sigma_0}^\sigma\arg F(x+it_0)\,dx-\int_{\sigma_0}^\sigma\arg F(x+iT)\,dx,
\end{aligned}
\end{equation}
where $\arg F(s)$ should be taken continuously in the broken line with vertices at $\sigma+it_0$, $\sigma_0+it_0$, $\sigma_0+iT$, $\sigma+iT$. 

By Lemma \ref{L:Jexpres} we have $\arg F(s)=\arg J(s)+\arg(1+U(s))$, with $s=\sigma_0+it$. By the election of $\sigma_0$, the function $\arg(1+U(\sigma_0+it))$  for $t_0\le t\le T$  is always between $-\pi/2$ and $\pi/2$. By \eqref{E:logJ} we easily see that   $\arg F(\sigma+it)$ for $1-\sigma\ge t^a$, $t\ge t_0$  is equal to 
\begin{multline*}\arg F(\sigma+it)=\frac{\pi\sigma}{4}+\frac{\pi}{2}-\arctan\frac{\sigma}{t}+\frac12\Bigl(\frac\pi2-\arctan\frac{1-\sigma}{t}\Bigr)-\frac{3\pi}{8}\\+\frac{\pi}{\sqrt{2\pi}}(t^2+(1-\sigma)^2)^{1/4}\cos(\tfrac12\arctan\tfrac{1-\sigma}{t})+\Orden(1).\end{multline*}
Since $1-\sigma_0=T^a$ it follows that $\arg F(\sigma_0+it)$ varies continuously between the extreme values
\[\arg F(\sigma_0+iT)=\pi\sqrt{\frac{T}{2\pi}}+\Orden(T^a)\qquad \arg F(\sigma_0+it_0)=\frac{\pi\sigma_0}{4}+\Orden(T^{a/2}).\]
Therefore, $\arg F(\sigma_0+it_0)=\Orden(T^a)$ and the lemma \ref{L:underside} apply to show that
\begin{equation}\label{intarc1}
\int_{\sigma_0}^{\sigma}\arg F(x+it_0)\,dx=\Orden(T^{2a}\log T).
\end{equation}
On the upper side for $s=x+iT$  with $\sigma_0\le x\le\sigma$ we have
\[\arg F(s)=\arg(s\pi^{-s/2}\Gamma(s/2))+\arg\Rzeta(s)\]
where the arguments on the right-hand side should be taken continuous and such that at the point $s_0=\sigma_0+iT$ we have 
\[
\arg(s_0\pi^{-s_0/2}\Gamma(s_0/2))+\arg\Rzeta(s_0)= \arg F(s_0)=\pi\sqrt{T/2\pi}+\Orden(T^a).
\]
We have some freedom here; we may pick $\arg(s_0\pi^{-s_0/2}\Gamma(s_0/2))=\Orden(T^a)$ so that Lemma \ref{L:upperside} applies, and 
\begin{equation}\label{intarc2}
\int_{\sigma_0+iT}^{\sigma+iT}\arg(s\pi^{-s/2}\Gamma(s/2)) \,ds=\Orden(T^{2a}).
\end{equation}
Then we must take $\arg\Rzeta(s_0)=\Orden(T^{1/2})$. We apply the Backlund lemma to bound the integral of $\arg\Rzeta(x+iT)$. We have 
\[|\arg\Rzeta(x+iT)-\arg\Rzeta(\sigma_0+iT)|=\Bigl|\Re\frac{1}{i}\int_{\sigma_0+iT}^{x+iT}\frac{\Rzeta'(z)}{\Rzeta(z)}\,dz\Bigr|\]

To determine the value of $\arg \Rzeta(x+iT)$, we apply Backlund's lemma.
Take a disc with center at $10+iT$ and radius $R=2(10-\sigma_0)$. On this disc
$|\Rzeta(s)|$ is bounded, according to Proposition 12 and 13 in \cite{A92} by 
\[\le C\max(T^{1/2}, T^{ cT^{a}}).\]
The value in the center of the disc $|\Rzeta(10+iT)|\ge 1/2$ by by  Proposition 6 in \cite{A173}. Hence, we get 
\[|\arg\Rzeta(x+iT)-\arg\Rzeta(\sigma_0+iT)|\le C T^a\log T.\]
It follows that 
\[|\arg\Rzeta(x+iT)|\le C T^a\log T+ CT^{1/2}\le CT^{1/2}.\]
From which we obtain 
\begin{equation}\label{intarc3}
\int_{\sigma_0+iT}^{\sigma+iT}\arg\Rzeta(s) \,ds=CT^{1/2+a}.
\end{equation}

Putting the results of  Propositions \ref{P:logF}, \ref{P:intvert} and equations \eqref{intarc1}, \eqref{intarc2} and \eqref{intarc3} into  equation \eqref{E:Litequation} yields 
\begin{align*}
2\pi\sum_{\substack{\beta\le\sigma\\t_0<\gamma\le T}}(\sigma-\beta)&=
-\frac{\pi T^2}{8}+(\sigma+1)\Bigl(\frac{T}{2}\log \frac{T}{2\pi}-\frac{T}{2}\Bigr)+\frac{T}{2}(\log 2+2\log(2\pi))-\frac{\pi\sigma^2}{8}\\&+\int_{t_0}^T\log|\Rzeta(\sigma+it)|\,dt+\orden(T^{1/2})+\frac{\pi T^2}{8}-\frac{T}{2}\log T+\frac{T}{2}-\frac{T}{2}\log2\pi\\ &+\Orden(T^{20/21})+\Orden(T^{2a}\log T)+\Orden(T^{\frac12+a}),
\end{align*}
Since $a=3/7$ the largest error term is $\Orden(T^{20/21})$. The term with $\sigma^2$ is less than this error term. Simplifying yields \eqref{E:thm2}.
\end{proof}

\begin{remark}
We can not apply Backlund's lemma in order to bound the integral of $\arg g(x+it_0)$ in Siegel's exposition \cite{Siegel}. To use $g(s)$ we  also have the extra difficulty that the Euler-MacLaurin approximation to the gamma function is not valid near the negative real axis. This is the main reason to use our $F(s)$ function instead of the  $g(s)$ function of Siegel.  
\end{remark}

\section{Proof of Propositions \ref{P:logF} and \ref{P:intvert}}\label{S:comp}

\subsection{Mean value of \texorpdfstring{$\log|\Gamma(s)|$}{loggamma}} We denote by $B_n(x)$ the Bernoulli polynomial and by $\widetilde{B_n}(x)$ the periodic function with period $1$ and such that $\widetilde{B_n}(x)=B_n(x)$ for $0< x<1$. They satisfy $B_n(x)=(-1)^n B_n(1-x)$. For an odd index greater than $1$ we have $B_{2n+1}(0)=B_{2n+1}(1/2)=B_{2n+1}(1)=0$, and these are the only zeros of $B_{2n+1}(x)$ on the interval $[0,1]$.

\begin{lemma}\label{L:simpleB}
Let $f\colon(a,b)\to[0,+\infty)$ be a positive and monotonous function, then 
\[\Bigl|\int_a^b f(x)\widetilde{B_3}(x)\,dx\Bigr|\le \frac{3}{64}\max(f(a), f(b)).\]
\end{lemma}
\begin{proof}
Assume that $f$ is not increasing, the proof in the other case is similar.
We have $B_3(x)=x(x-\frac12)(x-1)$, so $B_3(x)$ is positive in $(0,1/2)$ and negative in $(1/2,1)$. 
Define \[a_n=\left|\int_{n/2}^{(n+1)/2}f(x)\widetilde{B_3}(x)\,dx\right|=
(-1)^n\int_{n/2}^{(n+1)/2}f(x)\widetilde{B_3}(x)\,dx.\]
For any integer $k$, since $f(x+k)\ge f(1-x+k)$ for $0<x<\frac12$
\[a_{2k}=\int_{k}^{k+1/2}f(x)\widetilde{B_3}(x)\,dx=\int_0^{1/2}f(x+k)\widetilde{B_3}(x)\,dx\ge \int_0^{1/2}f(1-x+k)x(x-\tfrac12)(x-1)\,dx\]
changing variables $y=1-x$
\begin{multline*}=\int_{1/2}^1 f(y+k)(1-y)(\tfrac12-y)(-y)\,dy=-\int_{1/2}^1 f(y+k)\widetilde{B_3}(y)\,dy=\\
-\int_{k+1/2}^{k+1}f(y)\widetilde{B_3}(y)\,dy=a_{2k+1}.\end{multline*}
In the same way, we prove $a_{2k+1}\ge a_{2k+2}$. 

Suppose first that there exist $n <m$ such that  \[(n-1)/2<a\le n/2\le m/2\le b<(m+1)/2.\] If $n<m$  we will have
\[\int_a^b=\int_a^{n/2}+\int_{n/2}^{m/2}+\int_{m/2}^b=\int_a^{n/2}+\sum_{n\le k< m}(-1)^ka_k+\int_{m/2}^b ,\]
the omitted integrand being $f(x)\widetilde{B_3}(x)$, and with $a_n\ge a_{n+1}\ge\cdots\ge a_m\ge0$. Therefore,  we have 
\[\Bigl|\sum_{n\le k< m}(-1)^ka_k\Bigr|=a_n-a_{n+1}+a_{n+2}-\cdots \pm a_{m-1}\le a_n.\]
Therefore,
\[\Bigl|\int_a^b\Bigr|\le \Bigl|\int_a^{n/2}\Bigr|+a_n+\Bigl|\int_{m/2}^b\Bigr|\le 
3\int_{(n-1)/2}^{n/2} f(a) |\widetilde{B_3}(x)|\,dx=\frac{3 f(a)}{64}.\]
When $n=m$ or there is no fraction $n/2$ between $a$ and $b$, the inequality is easier to prove. If we assume that $f$ is increasing, we must substitute $f(a)$ by $f(b)$. 
So, in general, the inequality is true with the maximum between $f(a)$ and $f(b)$. 
\end{proof}
With the same procedure, we may prove this Proposition.
\begin{proposition}
Let $f\colon(a,b)\to[0,+\infty)$ be a positive and monotonous function, then for all $n\ge0$
\[\Bigl|\int_a^b f(x)\widetilde{B_{2n+1}}(x)\,dx\Bigr|\le (-1)^n3\Bigl(1-\frac{1}{2^{2n+2}}\Bigr)\frac{B_{2n+2}}{n+1}\max(f(a), f(b)).\]
\end{proposition}

\begin{proposition}
For $\sigma\in\R$ and $T>1$ we have 
\begin{equation}\label{E:meanloggamma}
\checked\int_1^T \log|\Gamma(\tfrac{\sigma+it}{2})|\,dt=\Im\int_{\sigma+i}^{\sigma+iT}\bigl((\tfrac{s-1}{2})\log \tfrac{s}{2}-\tfrac{s}{2}+\tfrac12\log2\pi+\tfrac{1}{6s}\bigr)\,ds+\Orden^*\Bigl(\frac{3\sqrt{3}}{16}\Bigr).
\end{equation}
\end{proposition}
\begin{proof}
Since $\sigma$ and $T>0$ are arbitrary, we cannot directly use the Euler-MacLaurin expansion. Instead, we use an intermediate expression obtained in the course of the proof of the Euler-MacLaurin expansion \cite{E}*{p.~109}. We have 
\[\log\Gamma(s)=(s-\tfrac12)\log s-s+\tfrac12\log2\pi+\frac{B_2}{2s}-\int_0^\infty\frac{\widetilde{B}_3(x)\,dx}{3(s+x)^3}, \qquad s+|s|\ne0.\]

Since
\[\int_1^T \log|\Gamma(\tfrac{\sigma+it}{2})|\,dt=\Im\int_{\sigma+i}^{\sigma+iT}
\log\Gamma(\tfrac{s}{2})\,ds,\]
we proceed to bound 
\[R:=\Im\int_{\sigma+i}^{\sigma+iT}\int_0^\infty\frac{\widetilde{B}_{3}(x)\,dx}{3(s/2+x)^3}\,ds.\]
The integral converge absolutely, so we may apply Fubini's Theorem. And we have 
\begin{align*}R&=\Im\int_0^\infty\widetilde{B}_{3}(x)\Bigl(\int_{\sigma+i}^{\sigma+iT}\frac{ds}{3(s/2+x)^3}\Bigr)\,dx\\&=\Im\frac43\int_0^\infty\widetilde{B}_{3}(x)\Bigl(\frac{1}{(\sigma+i+2x)^2}-\frac{1}{(\sigma+iT+2x)^2}\Bigr)\,dx\\
&=\frac43\int_0^\infty\widetilde{B}_{3}(x)\Bigl(-\frac{2(\sigma+2x)}{((\sigma+2x)^2+1)^2}+\frac{2T(\sigma+2x)}{((\sigma+2x)^2+T^2)^2}\Bigr)\,dx\end{align*}
Note that $\sigma\in\R$ may be negative. The function $\frac{2T y}{(y^2+T^2)^2}$
is negative and decreases in $(-\infty,-T/\sqrt{3})$, negative and increases in $(-T/\sqrt{3},0)$, positive and increasing in $(0,T/\sqrt{3})$ and decreasing and positive in $(T/\sqrt{3},+\infty)$
Its maximum value is achieved at $y=T/\sqrt{3}$ when it takes the value $\frac{3\sqrt{3}}{8T^2}$. At $-T/\sqrt{3}$ it has a minimum value of $-\frac{3\sqrt{3}}{8T^2}$. Hence, we may separate the integral into at most four integrals where Lemma \ref{L:simpleB} applies and 
\[\Bigl|\frac43\int_0^\infty\frac{2T(\sigma+2x)\widetilde{B}_{3}(x)}{((\sigma+2x)^2+T^2)^2}\,dx\Bigr|\le \frac43\cdot 4\cdot \frac{3}{64}\cdot\frac{3\sqrt{3}}{8T^2}=\frac{3\sqrt{3}}{32 T^2}.
\]
This is also valid for $T=1$, so we get $|R|\le \frac{3\sqrt{3}}{16}$. 

Then 
\[\int_1^T \log|\Gamma(\tfrac{\sigma+it}{2})|\,dt=\Im\int_{\sigma+i}^{\sigma+iT}\bigl((\tfrac{s-1}{2})\log \tfrac{s}{2}-\tfrac{s}{2}+\tfrac12\log2\pi+\tfrac{1}{6s}\bigr)\,ds+\Orden^*\Bigl(\frac{3\sqrt{3}}{16}\Bigr).\qedhere\]
\end{proof}
The same is true if we take the lower limit of the integrals equal to $t_0\ge1$. 

A computation of this integral gives us 
\begin{equation}\label{E:partialint}
\begin{aligned} 
\int_{t_0}^T &\log|\Gamma(\tfrac{\sigma+it}{2})|\,dt=
\Bigl(-\frac{T^2}{4}+\frac{\sigma^2}{4}-\frac{\sigma}{2}+\frac16\Bigr) \arg(\sigma+iT)\quad \checked\\&+\Bigl(\frac{\sigma T}{4}-\frac{T}{4}\Bigr)\log(\sigma^2+T^2)-\sigma\Bigl(\frac{3T}{4}+\frac{T\log2}{2}\Bigr)+
\frac{T}{2}+\frac{T}{2}\log2+\frac{T}{2}\log2\pi\\&+
\Bigl(\frac{t_0^2}{4}-\frac{\sigma^2}{4}+\frac{\sigma}{2}-\frac16\Bigr)\arg(\sigma+it_0)
+\Bigl(\frac{t_0}{4}-\frac{\sigma t_0}{4}\Bigr)\log(\sigma^2+t_0^2)+
\\&+\sigma\Bigl(\frac{3t_0}{4}+\frac{t_0\log2}{2}\Bigr)
-\frac{t_0}{2}-\frac{t_0}{2}\log2-\frac{t_0}{2}\log2\pi+
\Orden^*(3\sqrt{3}/16)
\end{aligned}
\end{equation}
where the argument always refers  to the main argument, since $t_0$ and $T>0$ for us $0<\arg(\sigma+it)<\pi$.

\begin{proposition}
Let $t_0\ge1$  and $0<a<\frac12$ be fixed real numbers. There is a positive  function $f\colon[0,+\infty)\to[0,+\infty)$ such that $f(T)=\orden(T^{1/2})$ and such that  given real numbers   $\sigma\le 10$ and $T> t_0$,  connected by $1-\sigma\le T^a$, then 
\begin{equation}\label{E:sigmaneg}
\begin{aligned}
\int_{t_0}^T \log|\Gamma(\tfrac{\sigma+it}{2})|\,dt&=
-\frac{\pi T^2}{8}+(\sigma-1)\Bigl(\frac{T}{2}\log \frac{T}{2}-\frac{T}{2}\Bigr)+\frac{T}{2}\log(2\pi)-\frac{\pi \sigma^2}{8}\\&+\Orden(f(T)).
\end{aligned}
\end{equation}
\end{proposition}
\begin{proof}
The value of the integral is given in \eqref{E:partialint}. 
Assume first that $-2t_0<\sigma\le 10$. In this case, many terms in \eqref{E:partialint} are bounded by a constant only depending on $t_0$. Eliminating these terms, we get 
\begin{align*} 
\int_{t_0}^T &\log|\Gamma(\tfrac{\sigma+it}{2})|\,dt=
-\frac{T^2}{4}\arg(\sigma+iT)+\Bigl(\frac{\sigma T}{4}-\frac{T}{4}\Bigr)\log(\sigma^2+T^2)\\&
-\sigma\Bigl(\frac{3T}{4}+\frac{T\log2}{2}\Bigr)+
\frac{T}{2}+\frac{T}{2}\log2+\frac{T}{2}\log2\pi+\Orden(t_0^2\log(1+t_0)).
\end{align*}
We have $\arg(\sigma+iT)=\frac{\pi}{2}-\arctan\frac{\sigma}{T}$ and 
$\log(\sigma^2+T^2)=2\log T+\log(1+\sigma^2/T^2)=2\log T+\Orden(T^{-2})$.
Expanding the $\arctan$, we get (with constants only depending on $t_0$)
\begin{align*}
\int_{t_0}^T \log|\Gamma(\tfrac{\sigma+it}{2})|\,dt&=-\frac{\pi T^2}{8}+\frac{T^2}{4}\Bigl(\frac{\sigma}{T}+\Orden(T^{-3})\Bigr)+\Bigl(\frac{\sigma T}{4}-\frac{T}{4}\Bigr)(2\log T+\Orden(T^{-2}))\\&
-\sigma\Bigl(\frac{3T}{4}+\frac{T\log2}{2}\Bigr)+
\frac{T}{2}+\frac{T}{2}\log2+\frac{T}{2}\log2\pi+\Orden(1).
\end{align*}
This is \eqref{E:sigmaneg} except for the term $-\sigma^2/8$, that in this case is contained in the error term. 

In the other case, when $\sigma_0\le\sigma\le -2t_0$, eliminating terms that are 
$\orden(T^{a})$ we get
\begin{align*} 
\int_{t_0}^T &\log|\Gamma(\tfrac{\sigma+it}{2})|\,dt=
\Bigl(-\frac{T^2}{4}+\frac{\sigma^2}{4}\Bigr) \arg(\sigma+iT)+\Bigl(\frac{\sigma T}{4}-\frac{T}{4}\Bigr)\log(\sigma^2+T^2)\\&-\sigma\Bigl(\frac{3T}{4}+\frac{T\log2}{2}\Bigr)+
\frac{T}{2}+\frac{T}{2}\log2+\frac{T}{2}\log2\pi
-\frac{\sigma^2}{4}\arg(\sigma+it_0)+\orden(T^{1/2})
\end{align*}

since $\sigma<0$ we have
\[\arg(\sigma+iT)=\frac{\pi}{2}-\arctan\frac{\sigma}{T},\quad \arg(\sigma+it_0)=\pi+\arctan\frac{t_0}{\sigma},\]\[\log(\sigma^2+T^2)=2\log T+\log(1+\sigma^2/T^2),\quad \log(\sigma^2+t_0^2)=2\log|\sigma|+\log(1+t_0^2/\sigma^2),\]
which we may expand in convergent power series since $\sigma/T$ and $t_0/\sigma$ are in absolute value $<1$. Removing terms less than $\orden(T^{1/2})$ yields \eqref{E:sigmaneg}.
\end{proof}
\subsection{Proof of Proposition \ref{P:logF}}

\begin{proof}[Proof of Proposition \ref{P:logF}]
By the definition of $F$ we have 
\[\int_{t_0}^T \log|F(\sigma+it)|\,dt=\int_{t_0}^T\log|(\sigma+it)\pi^{-\sigma/2}|\,dt+\int_{t_0}^T \log|\Gamma(\tfrac{\sigma+it}{2})|\,dt+\int_{t_0}^T \log|\Rzeta(\sigma+it)|\,dt\]
\begin{align*}\int_{t_0}^T\log|(\sigma+it)\pi^{-\sigma/2}|\,dt&=-\sigma\frac{T-t_0}{2}\log\pi+\frac{T}{2}\log(\sigma^2+T^2)-T+\sigma\arctan\frac{T}{\sigma}\\&-\frac{t_0}{2}\log(\sigma^2+t_0^2)+t_0-\sigma\arctan\frac{t_0}{\sigma}.\end{align*}
After some simplifications
\[\int_{t_0}^T\log|(\sigma+it)\pi^{-\sigma/2}|\,dt= T\log T-T-\sigma\frac{T}{2}\log\pi+\orden(T^{1/2}).\]
Combining this with \eqref{E:sigmaneg} yields
\begin{align*}
\int_{t_0}^T \log|F(\sigma+it)|\,dt&=
-\frac{\pi T^2}{8}+(\sigma+1)\Bigl(\frac{T}{2}\log \frac{T}{2\pi}-\frac{T}{2}\Bigr)+\frac{T}{2}(\log 2+2\log(2\pi))-\frac{\pi\sigma^2}{8}\\&+\int_{t_0}^T\log|\Rzeta(\sigma+it)|\,dt+\orden(T^{1/2}).\qedhere\end{align*}
\end{proof}

\subsection{Proof of Proposition \ref{P:intvert}}
Here we will assume that the exponent $a= 3/7$. The error term we will obtain is determined by the error in Theorem \ref{T:RzetaAprox} for $U(t)\ll t^{-1/21}$. This is not the best exponent, and its relation to $a=3/7$ is not as direct, as seen in \cite{A193}. Therefore, we will not be able to obtain in Proposition \ref{P:intvert} an error better than $\Orden(T^{20/21})$. We suspect that Proposition \ref{P:intvert} can be improved. Therefore, we will retain in our equations some terms less than the final error $\Orden(T^{20/21})$, always assuming that $a\le 3/7$. In this way, we will get a conjecture about further terms in Proposition \ref{P:sumbeta}.

\begin{definition}
We define a periodic function $\per\colon\R\to\R$ by 
\begin{equation}
\per(x)=\sum_{n=1}^\infty \frac{2\sin(2\pi  n x)}{n^2}-\sum_{n=1}^\infty (-1)^n\frac{\sin(4\pi n x)}{n^2}.
\end{equation}
\end{definition}

Recall that for a given $s$ we define $\eta=\sqrt{(s-1)/2\pi i}$, taking the root such that $\Re\eta+\Im\eta>0$.  
\begin{lemma}\label{L:per}
Let $t_0>1$ be a fixed real number. For $T\to+\infty$ we have 
\begin{equation}
\Im\int_{1+it_0}^{1+iT}\log\frac{1-e^{2\pi i\eta}}{1+e^{4\pi i\eta}}\,ds=-\sqrt{T/2\pi}\;\per(\sqrt{T/2\pi})+\Orden(1),
\end{equation}
where the implicit constant only depends on $t_0$. 
\end{lemma}
\begin{proof}
Here $s=1+it$ with $t_0<t<T$ and $\eta=\sqrt{t/2\pi}$, hence 
\[J:=\Im\int_{1+it_0}^{1+iT}\log\frac{1-e^{2\pi i\eta}}{1+e^{4\pi i\eta}}\,ds=
\Im i\int_{t_0}^{T}\log\frac{1-e^{2\pi i\eta}}{1+e^{4\pi i\eta}}\,dt=
\Re\int_{t_0}^{T}\log\frac{1-e^{2\pi i\eta}}{1+e^{4\pi i\eta}}\,dt\]
Change variables $t=2\pi\eta^2$, and let $\tau_0=\sqrt{t_0/2\pi}$ and $\tau=
\sqrt{T/2\pi}$, then 
\[J=\Re4\pi\int_{\tau_0}^\tau \eta\log\frac{1-e^{2\pi i\eta}}{1+e^{4\pi i\eta}}\,d\eta\]
Now we use the expansion
\begin{equation}\label{E:powerseries}
\log\frac{1-z}{1+z^2}=-\sum_{n=1}^\infty\frac{z^n}{n}+\sum_{n=1}^\infty(-1)^n\frac{z^{2n}}{n}.
\end{equation}
Then we have (it is not difficult to justify the integration term by term)
\begin{align*}
J&=4\pi\Re\Bigl\{-\sum_{n=1}^\infty\frac{1}{n}\int_{\tau_0}^\tau\eta e^{2\pi i n\eta}\,d\eta+\sum_{n=1}^\infty\frac{(-1)^n}{n}\int_{\tau_0}^\tau\eta e^{4\pi i n\eta}\,d\eta\Bigr\}\\
&=4\pi\Re\Bigl\{\Bigl.\sum_{n=1}^\infty\Bigl(-\frac{e^{2\pi in \eta}}{4\pi^2n^3}+\frac{i\eta e^{2\pi in \eta}}{2\pi n^2}+(-1)^n\frac{e^{4\pi in\eta}}{16\pi^2n^3}-(-1)^n\frac{i\eta e^{4\pi i n\eta}}{4\pi n^2}\Bigr)\Bigr|_{\tau_0}^\tau\Bigr\}
\end{align*}
All  terms in $\tau_0$ will contribute only to a term $\Orden(1)$ (with a constant dependent on $t_0$). Analogously, the term in $\tau$ and $n^3$,  which do not have the factor $\eta$ also contributes to a term $\Orden(1)$ (this time with absolute constants). So, we get 
\[J=\Re\sum_{n=1}^\infty\Bigl(\frac{2i\tau e^{2\pi i n\tau}}{n^2}-(-1)^n\frac{i\tau e^{4\pi i n\tau}}{n^2}\Bigr)+\Orden(1)=
\Re i\tau\sum_{n=1}^\infty\Bigl(\frac{2e^{2\pi i n\tau}}{n^2}-(-1)^n\frac{e^{4\pi i n\tau}}{n^2}\Bigr)+\Orden(1)\]
\[=\tau\sum_{n=1}^\infty\Bigl(-\frac{2\sin(2\pi n\tau)}{n^2}+(-1)^n\frac{\sin(4\pi n\tau)}{n^2}\Bigr)+\Orden(1)=-\tau\per(\tau)+\Orden(1).\qedhere\]
\end{proof}
The previous lemma considers the integral with limits $1+it_0$ and $1+iT$, but we are interested in this integral with limits $\sigma_0+it_0$ and $\sigma_0+iT$. The difference between the two is bounded in the next lemma.

\begin{lemma}\label{L:restmean}
Let $0<a\le 3/7$ and $1-\sigma\le T^a$, then for $T\to+\infty$ we have
\begin{equation}
\Im\int_{\sigma+it_0}^{1+it_0}\log\frac{1-e^{2\pi i\eta}}{1+e^{4\pi i\eta}}\,ds-\Im\int_{\sigma+iT}^{1+iT}\log\frac{1-e^{2\pi i\eta}}{1+e^{4\pi i\eta}}\,ds=\Orden(T^a).
\end{equation}
\end{lemma}
\begin{proof}
The two integrals are treated in the same way. For the second, for example, we have  $s=x+iT$, with $\sigma< x< 1$. Then as we noticed before $|e^{2\pi i \eta}|<1$ and the logarithm is defined by means of \eqref{E:powerseries} with $z=e^{2\pi i \eta}$. Both series in \eqref{E:powerseries} have imaginary parts in $(-\pi/2,\pi/2)$ so that 
\[
\Bigl|\Im\int_{\sigma+iT}^{1+iT}\log\frac{1-e^{2\pi i\eta}}{1+e^{4\pi i\eta}}\,ds\Bigr|=\Bigl|\int_{\sigma}^{1}\Im\log\frac{1-e^{2\pi i\eta}}{1+e^{4\pi i\eta}}\,dx\Bigr|\le \int_{\sigma}^{1}\pi\,dx=\pi(1-\sigma)=\Orden(T^a).\qedhere\]
\end{proof}

\begin{lemma}\label{L:littlemean}
Let $1\le t_0$ and $a\le \frac37$ be  fixed real numbers and  $1-\sigma_0= T^a$.

Then for $T\to+\infty$ we have
\begin{equation}
\int_{t_0}^T\log|e^{\pi i\eta}|\,dt=-\pi(1-\sigma_0)\Bigl(\frac{T}{2\pi}\Bigr)^{1/2}+\frac{\sqrt{\pi}}{3}(1-\sigma_0)^{3/2}
+\Orden(T^{a/2}).
\end{equation}
\end{lemma}
\begin{proof}
Put $s=\sigma_0+it$, then 
\[\int_{t_0}^T\log|e^{\pi i\eta}|\,dt=\Im\int_{\sigma+it_0}^{\sigma+iT}\log e^{\pi i \eta}\,ds=\Im\int_{\sigma+it_0}^{\sigma+iT}\pi i\sqrt{\frac{s-1}{2\pi i}}\,ds\]
\begin{align*}
&=-\Im \int_{\sigma+it_0}^{\sigma+iT} 2\pi^2\sqrt{\frac{s-1}{2\pi i}}\,d\frac{s-1}{2\pi i} =-2\pi^2\Im \Bigl.\frac{1}{3/2}\Bigl(\frac{s-1}{2\pi i}\Bigr)^{3/2}\Bigr|_{s=\sigma+it_0}^{\sigma+iT}\\
&=-\frac{4\pi^2}{3}\Im\Bigl\{\Bigl(\frac{T}{2\pi}+i\frac{1-\sigma}{2\pi}\Bigr)^{3/2}-\Bigl(\frac{t_0}{2\pi}+i\frac{1-\sigma}{2\pi}\Bigr)^{3/2}\Bigr\}\\
&=\Im\Bigl\{-\frac{4\pi^2}{3}\Bigl(\frac{T}{2\pi}\Bigr)^{3/2}\Bigl(1+i\frac{1-\sigma}{T}\Bigr)^{3/2}+\frac{4\pi^2}{3}\Bigl(i\frac{1-\sigma}{2\pi}\Bigr)^{3/2}\Bigl(1-i\frac{t_0}{1-\sigma}\Bigr)^{3/2}\Bigr\}.
\end{align*}
According to the hypothesis, both $(1-\sigma)/T$ and $t_0/(1-\sigma)$ are less than $\frac12$, say. So, we may expand both parenthesis in convergent power series.
If I want to get this with an error $\orden(1)$ we must retain the terms  until $((1-\sigma)/T)^3$ and $(t_0/(1-\sigma))^2$.
Therefore, the result is the imaginary part of 
\begin{align*}
\Im\Bigl\{-\frac{4\pi^2}{3}\Bigl(\frac{T}{2\pi}\Bigr)^{3/2}\Bigl(1+\frac{3i}{2}\frac{1-\sigma}{T}-\frac{3}{8}\frac{(1-\sigma)^2}{T^2}\Bigr)+\Orden((1-\sigma)^3T^{-3/2})\\
+(-1+i)\frac{4\pi^2}{3\sqrt{2}}\Bigl(\frac{1-\sigma}{2\pi}\Bigr)^{3/2}
\Bigl(1-\frac{3i}{2}\frac{t_0}{1-\sigma}\Bigr)+\Orden((1-\sigma)^{-1/2})\Bigr\}
\end{align*}
\begin{multline*}\int_{t_0}^T\log|e^{\pi i\eta}|\,dt=-\frac{4\pi^2}{3}\Bigl(\frac{T}{2\pi}\Bigr)^{3/2}\frac{3}{2}\frac{1-\sigma}{T}+\frac{4\pi^2}{3\sqrt{2}}\Bigl(\frac{1-\sigma}{2\pi}\Bigr)^{3/2}\\+\frac{4\pi^2}{3\sqrt{2}}\frac{3}{2}\Bigl(\frac{1-\sigma}{2\pi}\Bigr)^{3/2}\frac{t_0}{1-\sigma}+\Orden((1-\sigma)^3T^{-3/2})+\Orden((1-\sigma)^{-1/2})\end{multline*}
\[=-\pi(1-\sigma)\Bigl(\frac{T}{2\pi}\Bigr)^{1/2}+\frac{\sqrt{\pi}}{3}(1-\sigma)^{3/2}
+\frac{t_0}{2}(\pi(1-\sigma))^{1/2}+\Orden((1-\sigma)^{-1/2})\]
assuming that $a\le \frac37$ so that the exponent $3a-3/2\le-a/2$ . 
\end{proof}

\begin{proposition}\label{P:meanJsigma0}
Let $t_0>1$ greater than the constant in Theorem \ref{T:RzetaAprox} and $1-\sigma_0=T^a$ with $0<a<1/2$. Then, for $T\to+\infty$
\begin{equation}
\int_{t_0}^T\log|J(\sigma_0+it)|\,dt=-\frac{\pi T^2}{8}+\frac{T}{2}\log T-\frac{T}{2}+\frac{T}{2}\log (2\pi)+\Orden(T^{a+\frac12})
\end{equation}
\end{proposition}
\begin{proof}
By Lemmas \ref{L:Jexpres}, \ref{L:per}, \ref{L:restmean} and \ref{L:littlemean}, it remains to compute 
\[A=\Im\int_{\sigma_0+it_0}^{\sigma_0+iT}\Bigl(\frac{\pi i s}{4}+\log s-\frac{1}{2}\log\frac{1-s}{2\pi }-\frac{3\pi i}{8}\Bigr)\,ds\]
\[A=\Im\Bigl\{\Bigl.\frac{\pi i s^2}{8}+s \log s-\frac12s-\frac{3\pi i s}{8}+\frac{1-s}{2}\log\frac{1-s}{2\pi}\Bigr|_{\sigma_0+it_0}^{\sigma_0+iT}\Bigr\}\]
\begin{multline*}=-\frac{\pi T^2}{8}-\frac{T}{2}+\Orden_{t_0}(1)+\Im\Bigl\{(\sigma_0+iT)\log(\sigma_0+i T)-(\sigma_0+it_0)\log(\sigma_0+i t_0)\\+\frac{1-\sigma_0-iT}{2} \log\frac{1-\sigma_0-iT}{2\pi}-\frac{1-\sigma_0-it_0}{2} \log\frac{1-\sigma_0-i t_0}{2\pi}\Bigr\}\end{multline*} 
Simple computations yields (I have written here $\sigma$ instead of $\sigma_0$)
\[\text{\checked}\quad\Im(\sigma+iT)\log(\sigma+iT)=T\log T-T\sum_{n=1}^\infty\frac{(-1)^n}{2n}\Bigl(\frac{\sigma}{T}\Bigr)^{2n}+\frac{\pi\sigma}{2}+\sigma \sum_{n=1}^\infty \frac{(-1)^n}{2n-1}\Bigl(\frac{\sigma}{T}\Bigr)^{2n-1}.\]
\[\text{\checked}\quad\Im (\sigma+it_0)\log(\sigma+it_0)=t_0\log|\sigma|-t_0\sum_{n=1}^\infty \frac{(-1)^n}{2n}\Bigl(\frac{t_0}{\sigma}\Bigr)^{2n}+\pi \sigma-\sigma\sum_{n=1}^\infty \frac{(-1)^n}{2n-1}\Bigl(\frac{t_0}{\sigma}\Bigr)^{2n-1}.
\]
\begin{multline*}
\text{\checked}\quad\Im\Bigl(\frac{1-\sigma}{2}-\frac{iT}{2}\Bigr)\log\frac{1-\sigma-iT}{2\pi}\\=-\frac{T}{2}\log\frac{T}{2\pi}+\frac{T}{2}\sum_{n=1}^\infty\frac{(-1)^n}{2n}\Bigl(\frac{1-\sigma}{T}\Bigr)^{2n}-\frac{\pi(1-\sigma)}{4}-\frac {1-\sigma}{2}\sum_{n=1}^\infty\frac{(-1)^n}{2n-1}\Bigl(\frac{1-\sigma}{T}\Bigr)^{2n-1}.\end{multline*}

\begin{multline*}\hbox{\checked}\quad
\Im\Bigl(\frac{1-\sigma}{2}-\frac{it_0}{2}\Bigr)\log\frac{1-\sigma-it_0}{2\pi}\\=-\frac{t_0}{2}\log\frac{1-\sigma}{2\pi}+\frac{t_0}{2}\sum_{n=1}^\infty\frac{(-1)^n}{2n}\Bigl(\frac{t_0}{1-\sigma}\Bigr)^{2n}+\frac{1-\sigma}{2}
\sum_{n=1}^\infty
\frac{(-1)^n}{2n-1}\Bigl(\frac{t_0}{1-\sigma}\Bigr)^{2n-1}\Bigr).\end{multline*}
Substituting these values, we obtain 
\[A=-\frac{\pi T^2}{8}-\frac{T}{2}+T\log T-\frac{T}{2}\log\frac{T}{2\pi}+\Orden(T^a).\]
\[A=-\frac{\pi T^2}{8}+\frac{T}{2}\log T-\frac{T}{2}+\frac{T}{2}\log2\pi+\Orden(T^a).\]
By the above lemmas \ref{L:Jexpres}, \ref{L:per}, \ref{L:restmean} and \ref{L:littlemean} we get then 
\begin{align*}
\int_{t_0}^T&\log|J(\sigma_0+it)|\,dt=-\frac{\pi T^2}{8}+\frac{T}{2}\log T-\frac{T}{2}+\frac{T}{2}\log2\pi+\Orden(T^a)\\
&-\pi(1-\sigma_0)\Bigl(\frac{T}{2\pi}\Bigr)^{1/2}+\frac{\sqrt{\pi}}{3}(1-\sigma_0)^{3/2}+\Orden(T^{a/2})\\
&-\Bigl(\frac{T}{2\pi}\Bigr)^{1/2}\per(\sqrt{T/2\pi})+\Orden(1).
\end{align*}
Therefore, 
\begin{equation}\label{E:extraterms}
\begin{aligned}
\int_{t_0}^T\log|J(\sigma_0+it)&|\,dt=-\frac{\pi T^2}{8}+\frac{T}{2}\log T-\frac{T}{2}+\frac{T}{2}\log2\pi+\pi \sigma_0\Bigl(\frac{T}{2\pi}\Bigr)^{1/2}\\&-\bigl(\pi+\per(\sqrt{T/2\pi})\bigr)\Bigl(\frac{T}{2\pi}\Bigr)^{1/2}+\frac{\sqrt{\pi}}{3}(1-\sigma_0)^{3/2}+\Orden(T^a).
\end{aligned}
\end{equation}
\end{proof}
\begin{remark}
Siegel \cite{Siegel}*{p.~304} asserts that he can prove his theorems with $a=\varepsilon$. This would imply that the term containing the function $P$ in equation \eqref{E:extraterms} would make sense. In \cite{A172} we saw that this term appears in the experimental data obtained with the zeros of $\Rzeta(s)$.
\end{remark}

\begin{proof}[Proof of Proposition \ref{P:intvert}]  Since $t_0$ is chosen adequately, we have by Theorem \ref{T:RzetaAprox} that $F(s)=J(s)(1+U(s))$, with $|U(s)|\le At^{-1/21}$.  It follows that 
\[\int_{t_0}^T\log|1+U(\sigma_0+it)|\,dt=\Orden(T^{20/21}).\]
Combining this with Proposition \ref{P:meanJsigma0} proves Proposition \ref{P:intvert}.
\end{proof}

\end{document}